\documentclass[11pt,a4paper]{amsart} 
\usepackage{url} 
\usepackage{mathrsfs}
\usepackage{enumitem}
\usepackage{color} 
\usepackage[colorlinks,linkcolor=blue,anchorcolor=blue,citecolor=blue,backref=page]{hyperref}
\usepackage{graphics,epsfig}
\usepackage{graphicx}
\usepackage{float}
\usepackage{epstopdf}
\usepackage[utf8]{inputenc}
\usepackage{hyperref} 
\usepackage{calc, amscd}
\usepackage[T1]{fontenc} 
\usepackage{pdfpages}
\usepackage{mathtools} 
\usepackage{scalerel,stackengine}
\usepackage{caption}
\usepackage{mathptmx}
\usepackage[per-mode=symbol]{siunitx}
\usepackage{bm}
\usepackage{amsmath}
\hypersetup{breaklinks=true}
\usepackage{amsthm,mathrsfs,multirow,xcolor,lscape,longtable,pbox,lipsum}
\usepackage{amsfonts,amscd,amssymb}
\usepackage[thinlines]{easytable}
\usepackage{microtype}
\usepackage{todonotes}
\usepackage{phaistos}
\usepackage{tabularx}
\newtheorem{theorem}{Theorem}

\newtheorem{lemma}[theorem]{Lemma}

\newtheorem{proposition}{Proposition}

\newtheorem*{theorem*}{Theorem}
\newtheorem*{question*}{Question}

\usepackage{mathtools}

\def\exp{\mathrm{exp}}

\begin{document}
\date{\today}
\author{R. Muneeswaran, S. Krishnamoorthy and Subham Bhakta  }

\address{Indian Institute of Science Education and Research, Thiruvananthapuram, India} 
\email{muneeswaran20@iisertvm.ac.in}
\address{Indian Institute of Science Education and Research, Thiruvananthapuram, India} 
\email{srilakshmi@iisertvm.ac.in}
\address{School of Mathematics \& Statistics, University of New South Wales, Australia.} 
\email{subham.bhakta@unsw.edu.au}
\newcommand\G{\mathbb{G}}
\newcommand\sO{\mathcal{O}}
\newcommand\sE{{\mathcal{E}}}
\newcommand\tE{{\mathbb{E}}}
\newcommand\sF{{\mathcal{F}}}
\newcommand\sG{{\mathcal{G}}}
\newcommand\sH{{\mathcal{H}}}
\newcommand\sN{{\mathcal{N}}}
\newcommand\GL{{\mathrm{GL}}}
\newcommand\HH{{\mathrm H}}
\newcommand\mM{{\mathrm M}}
\newcommand\fS{\mathfrak{S}}
\newcommand\fP{\mathfrak{P}}
\newcommand\fQ{\mathfrak{Q}}
\newcommand\Qbar{{\bar{\Q}}}
\newcommand\sQ{{\mathcal{Q}}}
\newcommand\sP{{\mathbb{P}}}
\newcommand{\Q}{\mathbb{Q}}
\newcommand{\tH}{\mathbb{H}}
\newcommand{\Z}{\mathbb{Z}}
\newcommand{\R}{\mathbb{R}}
\newcommand\Gal{{\mathrm {Gal}}}
\newcommand\SL{{\mathrm {SL}}}
\newcommand\Hom{{\mathrm {Hom}}}

\def \bphi{\bm{\varphi}}

\newtheorem{thm}{Theorem}[section]
\newtheorem{ack}[thm]{Acknowledgement}

\newtheorem{cor}[thm]{Corollary}
\newtheorem{conj}{Conjecture}
\newtheorem{prop}[thm]{Proposition}
\newtheorem{quest}{Question}

\theoremstyle{definition}

\newtheorem{claim}[thm]{Claim}

\theoremstyle{remark}

\newtheorem*{fact}{Fact}
\newcommand{\SK}[1]{{\color{red} \sf Srilakshmi: [#1]}}
\newcommand{\RM}[1]{{\color{blue} \sf  Muneeswaran: [#1]}}
\newcommand{\SB}[1]{{\color{blue} \sf  Subham: [#1]}}

\title[Indivisibility in quadratic number fields]{ON AN INDIVISIBILITY VERSION OF IIZUKA'S CONJECTURE}
\begin{abstract}
Iizuka's conjecture predicts that, given $m \in \mathbb{N}$ and a prime $p$, there exists infinitely many integers $n$ such that the class numbers of \textit{all} of the following quadratic number fields,
\[
\mathbb{Q}(\sqrt{n}),\ \mathbb{Q}(\sqrt{n+1}),\ \ldots,\ \mathbb{Q}(\sqrt{n+m}),
\]
are divisible by $p$. In this article, given $k$ and $m$, we study the proportion of $n$ such that the class numbers of \textit{none} of the successive fields
\[
\mathbb{Q}(\sqrt{n}),\ \mathbb{Q}(\sqrt{n+1}),\ \ldots,\ \mathbb{Q}(\sqrt{n+m}),
\]
are divisible by \( 3^k \). Moreover, we study the proportion of imaginary biquadratic fields whose class numbers are not divisible by $3$.
\end{abstract}

\maketitle
\tableofcontents

\section{Introduction}
\subsection{Set-up and motivation}
Let $K$ be a number field and let $\mathcal{O}_K$ denote its ring of integers. The class group of $K$ is a finite abelian group, and its order is known as the class number $h_K$ of $K$, which plays an important role in understanding the arithmetic of $K$. It is known that the only imaginary quadratic fields with class numbers exactly equal to one are $\Q(\sqrt{-d}), \ d=1,2,3,7,11,19,43,67,163$. Gauss conjectured that there are infinitely many real quadratic fields with class numbers exactly equal to one.

Examining the divisibility properties of class numbers is essential for understanding the structure of class groups. For a fixed $n_0$, many authors have identified infinitely many quadratic fields with class numbers divisible by $n_0$ \cite{MR85301,KriPasMun,MR3069394,MR266898,MR335471,MR4270672,louboutin2009divisibility,krishnamoorthy2021note,hoque2017divisibility,murty1999exponents,MR2843095,MR4761048}.  In fact, the Cohen-Lenstra heuristics \cite{cohen1984heuristics} predicts that the proportion of quadratic fields with class numbers divisible by $n_0$, is positive. Iizuka \cite{MR3724158} produced infinitely many $n$ such that the class numbers of $\Q(\sqrt{n})$ and $ \Q(\sqrt{n+1})$ are divisible by $3$. Moreover, he conjectured the following. 

\begin{conj}{(Iizuka)}
Given $m\in\mathbb{N}$ and a prime $p$, there exist infinitely many integers $n$ such that the class numbers of all of the following quadratic number fields,
$$\Q(\sqrt{n}),\ \Q(\sqrt{n+1})\ldots,\ \Q(\sqrt{n+m})$$ are divisible by $p$.
\end{conj}

The main goal of this article is to study the \textit{indivisibility} properties of class numbers. Gauss proved that class number of $\Q(\sqrt{-p}), \ p\equiv3\pmod{4}$ is odd. This gives us infinitely many quadratic fields whose class numbers are not divisible by $2$. Hartung \cite{MR352040} proved the existence of an infinite family of imaginary quadratic fields whose class numbers are not divisible by $3$. 

One of the significant aspects of the indivisibility of class numbers by a prime $p$ is understanding the Iwasawa invariants $\lambda_p(K)$ and $\mu_p(K)$ associated with the cyclotomic $\Z_p$ extension over $K$. Ferrero and Washington~\cite{FW} showed that $\mu_p(K)=0$ for any abelian extension $K/\Q$. It is well-known by Iwasawa~\cite{iwasawa} that, for any quadratic $K/\Q$, if $p$ does not split completely in $K$ and $p$ does not divide the class number $h_K$, then $\lambda_p(K)=0$. Greenberg~\cite{taya} conjectured that $\lambda_p(K)=0$ for any totally real extension $K/\Q$. This article produces several quadratic extensions of $\Q$ for which the $\lambda_3$ invariant is zero.

Byeon~\cite{MR2073286} proved that for a fixed square-free integer $t$, there exists infinitely many $d$ such that the class numbers of both $\Q(\sqrt{d})$ and $ 
 \Q(\sqrt{td})$ are not divisible by $3$. One of the interesting applications of Byeon's results is in extending Vatsal~\cite{vatsal}, by providing a positive proportion of rank 1 twists for $X_0(19)$.

\subsection{Main results}
Motivated by Iizuka’s conjecture and by results on the indivisibility of class numbers, we study the problem of \textit{simultaneous indivisibility} of class numbers of quadratic number fields. In this direction, Chattopadhyay and Saikia \cite[Theorem~1]{CS22} showed that there exists a set of positive density 
consisting of positive fundamental discriminants $d$ such that the class numbers of both 
$\Q(\sqrt{d})$ and $\Q(\sqrt{d+4})$ are not divisible by $3$. 

A natural extension of their result is to ask whether a similar indivisibility condition 
can hold simultaneously for the quadratic fields associated with five consecutive integers. 
Since in this case not all of these integers correspond to fundamental discriminants, 
we instead consider the problem of finding natural numbers $n$ such that 
\textit{none} of the class numbers of 
\[
\Q(\sqrt{n}),\ \Q(\sqrt{n+1}),\ \ldots,\ \Q(\sqrt{n+4})
\] 
are divisible by $3$.

In fact, below, we establish a more general result, in which we study the indivisibility of class numbers by any power of $3$, and also with a higher number of consecutive terms.

\begin{theorem}\label{indizuka}
Let $k>0$ and $m\geq 0$ be two fixed integers. Then the set of natural numbers \( n \in \mathbb{N} \) such that none of the class numbers of the real quadratic fields
\[
\mathbb{Q}(\sqrt{n}), \, \mathbb{Q}(\sqrt{n+1}), \, \ldots,\ \mathbb{Q}(\sqrt{n+m}),
\]
are divisible by \( 3^k \), has density at least $1-\frac{m+1}{3(3^k-1)}$ in $\mathbb{N}$.
\end{theorem}

In particular, Theorem~\ref{indizuka} addresses the question raised earlier for $k=1$. For any fixed \( m \), as one should perhaps expect, Theorem~\ref{indizuka} shows that the density of \( m \)-tuples of consecutive real quadratic fields whose class numbers are indivisible by \( 3^k \), tends to \(1\) as \(k \to \infty\). Moreover, \cite[Conjecture]{BM19} predicts the precise densities at least for $m=0$, and in particular asserts that while the density increases strictly with $k$, it always remains strictly less than $1$. Therefore, Theorem~\ref{indizuka} provides some meaningful lower bounds in the light of \cite[Conjecture]{BM19}; in the sense that, for any fixed \(k_0\), quadratic fields whose class numbers are not divisible by \(3^{k_0}\), do not necessarily serve as examples of quadratic fields whose class numbers are indivisible by \(3^k\) for any \(k > k_0\).

For a quick example, note that none of the class numbers of  
\[
\Q(\sqrt{443}),\ \Q(\sqrt{444}),\ \Q(\sqrt{445}),\ \dots,\ \Q(\sqrt{464})
\]  
is divisible by $3^2$. Thus, for $k=2$, we may take $n=443$ and $m=21$. 
However, since the class number of $\Q(\sqrt{443})$ is $3$, we cannot take $n=443$ when $k=1$.

Now, for any integer $i\geq 1$, let us denote $d_i$ to be the $i^{th}$ square-free natural number. Then, we have the following analog of Theorem~\ref{indizuka}.
\begin{theorem}\label{squrefree}
Let $k>0$ and $m\geq 0$ be two fixed integers. Then the set of natural numbers \( i \in \mathbb{N} \) such that none of the class numbers of the real quadratic fields fields $$\Q(\sqrt{d_i}), \ \Q(\sqrt{d_{i+1}}),\ldots,\ \Q(\sqrt{d_{i+m}}),$$ 
are divisible by \( 3^k \), has density at least $1-\frac{m+1}{3(3^k-1)}$ in $\mathbb{N}$.
\end{theorem}

The proofs of both Theorem~\ref{indizuka} and Theorem~\ref{squrefree} are essentially consequences of an estimate for the average $3$-torsion parts of quadratic fields (see Lemma~\ref{prop} and Lemma~\ref{avergesum in N}) over suitable sets of integers. For this, we obtain an explicit version of Nakagawa and Horie~\cite{MR0958035} in Proposition~\ref{error}, using Bebalas~\cite[Theorem 2.3]{Belabas}.

Of course, both Theorem~\ref{indizuka} and Theorem~\ref{squrefree} are about real quadratic fields. To build a bridge with the imaginary quadratic fields, we use Scholz~\cite{MR1581309} (see Theorem~\ref{scholz}), and study (in Proposition~\ref{prop:alldivby3}) the average of $3$-torsions parts of the quadratic fields $\Q(\sqrt{3n})$. Consequently, we deduce the following result.
\begin{theorem}\label{indizukaimaginary}
Let $k>0$ and $m\geq 0$ be two fixed integers. Then the set of natural numbers \( n \in \mathbb{N} \) such that none of the class numbers of the real quadratic fields
\[
\mathbb{Q}(\sqrt{-n}), \, \mathbb{Q}(\sqrt{-n+1}), \, \ldots,\ \mathbb{Q}(\sqrt{-n+m}),
\]
are divisible by \( 3^{k+1} \), has density at least $1-\frac{m+1}{3(3^k-1)}$ in $\mathbb{N}$.
\end{theorem}

Certainly, Theorem~\ref{indizukaimaginary} is a weaker analog of Theorem~\ref{indizuka}, as Theorem~\ref{indizukaimaginary} is about indivisibility by a higher power of $3$. The limitations arise, as we do not explicitly know when we have $r=s$ and when $r=s+1$, in Theorem~\ref{scholz}.

Furthermore, we study the indivisibility problem for imaginary bi-quadratic fields. Note that any such field can be written as $\Q(\sqrt{-n_1},\sqrt{n_2})$, where $(n_1,n_2)\in-\mathbb{N}\times\mathbb{N}$. 

We denote the following counting function
\begin{equation}\label{eq:H}
H(X,Y)= \sharp\, \left\{ (n_1, n_2) \in \mathbb{N} \times \mathbb{N} : n_1 \leq X, \ n_2 \leq Y,~3 \nmid h_{\mathbb{Q}(\sqrt{-n_1}, \sqrt{n_2})} \  \right\},
\end{equation}
and prove the following estimate,

\begin{theorem}\label{biquad}
There exists a constant \( C > 0 \) such that
\[
H(X, Y) \geq CXY
\]
for any sufficiently large \( X, Y \) satisfying
\[
\log X = o\left( \sqrt{ \log Y \cdot \log \log Y } \right).
\]
\end{theorem}

Note that the condition on $X,Y$ holds, for instance when $Y\ge \exp((\log X)^c)$ for some $c>2$. In other words, Theorem~\ref{biquad} produces a positive proportion of imaginary bi-quadratic fields in \textit{thin rectangles}, whose class numbers are not divisible by $3$. 

\subsection{Notations}
As usual, by $O_{S_1, S_2,\ldots, S_k}(B)$ we mean a quantity with absolute value at most $cB$ for some positive constant $c$ depending on $S_1,\cdots, S_K$ only; if the subscripts are omitted, the implied constant is absolute. Also, we write $A=o(B)$ for $A/B\to 0$.

For a given set $\mathbf{S}$, we denote $\sharp\, \mathbf{S}$ as its cardinality. Moreover, we use the standard notations $\varphi(n)$ for the Euler function, $\mu(n)$ for the M{\"o}bius functions, and $(a,b)$ to denote the greatest common divisor of integers $a$ and $b$. Also, as usual, we use $\nu_p(n)$ to denote the $p$-adic valuation of a number $n\in \mathbb{N}$.

For any integer $n$ and prime $p$, we denote $r_p(n)$ to be the $p$-rank of $\Q(\sqrt{n})$. Also, we denote $S$ to be the set of all squarefree natural numbers, and $S_{+}$ (resp. $S_{-}$) to be the set of all positive (resp. negative) fundamental discriminants.

\section{Preliminaries}
In this section, we collect the preliminary results that are needed to prove the main results of this article.

Let us recall that an integer $d$ is called a \textit{fundamental discriminant}, if $d$ does not have an odd square factor other than $1$, and 
$$d \equiv 1 \pmod 4, \quad \mathrm{or}\quad d \equiv 8~\mathrm{or}~12 \pmod{16}.$$

In particular, let us note the following simple fact, which will be used many times in this article
\begin{equation}\label{eq:simplefact}
    \text{If}\quad d \in S \quad \mathrm{and}\quad d \not\equiv 1 \pmod{4}, \quad \text{ then}\quad 4d \in S_{+}.
\end{equation}

\subsection{Cohen-Lenstra's heuristics}
Let $p$ be any odd prime. Cohen-Lenstra~\cite[Page~57]{cohen1984heuristics} (see also \cite[Conjecture~4]{BM19}) conjectured the following.

\begin{conj}\label{cohen}
Let $S_+(X)$ (resp.\ $S_-(X)$) be the set of all positive (resp.\ negative) fundamental discriminants with $|d| \le X$. Then  
\[
\lim_{X\to\infty} \frac{\sum\limits_{d\in S_+(X)} p^{r_p(d)}}{\sum\limits_{d\in S_+(X)} 1} 
= 1 + \frac{1}{p}, 
\quad\text{and}\quad 
\lim_{X\to\infty} \frac{\sum\limits_{d\in S_-(X)} p^{r_p(d)}}{\sum\limits_{d\in S_-(X)} 1} 
= 2.
\]
\end{conj}

The conjecture is proved for $p=3$ by Davenport and Heilbronn~\cite{MR0491593}. We also note an extension by Nakagawa and Horie~\cite{MR0958035} over congruence classes, which will be a key ingredient in the proof of Theorem \ref{indizuka}.
\begin{lemma}[Nakagawa and
Horie~\cite{MR0958035}]\label{Horie}
Let $m$ and $N$ be two natural numbers, satisfying the following conditions:
    \begin{itemize}
        \item[(i)] If $p$ is an odd prime such that $p\mid\gcd(m,N)$, then $N\equiv0\pmod{p^2}$ and $m\not\equiv0\pmod{p^2}$.
        \item[(ii)] If $N$ is an even integer, then either $N\equiv0\pmod{4}$ and $m\equiv1\pmod{4}$, or $N\equiv0\pmod{16}$ and $m\equiv8,12\pmod{16}$.
    \end{itemize}
     Let us denote 
     $$S_+(X,m,N)=\{d\in S_+(X): d\equiv m\pmod{N}\},\quad \mathrm{and}$$ 
     $$S_-(X,m,N)=\{d\in S_-(X): d\equiv m\pmod{N}\}.$$ 
     Then we have,
$$\lim\limits_{X\to\infty}\frac{\sum\limits_{d\in S_+(X,m,N)}3^{r_3(d)}}{\sharp\, S_+(X,m,N)}=\frac{4}{3},\quad \mathrm{and}\quad \lim\limits_{X\to\infty}\frac{\sum\limits_{d\in S_-(X,m,N)}3^{r_3(d)}}{\sharp\, S_-(X,m,N)}=2.$$  
\end{lemma}

Using Lemma~\ref{Horie}, Nakagawa and
Horie deduced the following estimates.
\begin{lemma}[Nakagawa and
Horie~\cite{MR0958035}]\label{densS+}
For any $m,N$ as in Lemma~\ref{Horie}, we have
 $$   \liminf\limits_{X\to\infty}\frac{\sharp\, \left\{d\in S_+(X,m,N): 3\nmid h_{\Q(\sqrt{d})}\right\}}{\sharp\, S_+(X,m,N)}\geq \frac{5}{6},\quad \mathrm{and}$$
 $$ \liminf\limits_{X\to\infty}\frac{\sharp\, \left\{d\in S_-(X,m,N): 3\nmid h_{\Q(\sqrt{d})}\right\}}{\sharp\,  S_-(X,m,N)}\geq \frac{1}{2}.$$
\end{lemma}

Note that $m=1$ and $N=1$ satisfy the conditions in the Lemma \ref{Horie}, and in this case we have $S_+(X,m,N)=S_+(X)$.

\subsection{An explicit version of Nakagawa and Horie} The main goal of this section is to make Lemma~\ref{Horie} explicit. To begin with, we present an asymptotic formula for \( S_+(X, m, N) \) in the following lemma. While this result may exist in the literature, we include a proof here due to the absence of a clear reference and for the sake of completeness.
\begin{lemma}\label{S(X) error}
Let $m$ and $N$ be as in Lemma~\ref{Horie}. Then for any $\varepsilon>0$, we have the following estimate 
\begin{align*}
\sharp\, S_+(X,m,N) = \frac{3X}{\pi^2\varphi(N)}\prod\limits_{p\mid N}\frac{q}{1+p} + O_{\varepsilon}\left(X^{1/2}+N^{1/2+\varepsilon}\right),
\end{align*}
$\mathrm{where}$\quad 
\[
q = \begin{cases}
4 & \text{if } p=2 \\
p & \text{otherwise}.
\end{cases}
\] 
\end{lemma}
\begin{proof}
For any two natural numbers $a,q$ with $(a,q)=g\geq 1$, note the following well known estimate (see \cite[Equation (1.3)]{Nunes}, for instance)
\begin{align}\sum\limits_{\substack{1\leq n\leq X\\ n\equiv a\pmod{q}}}\mu^2(n) &= \sum\limits_{\substack{1\leq n\leq \frac{X}{g}\\ n\equiv \frac{a}{g}\pmod{\frac{q}{g}}}}\mu^2(n)= \frac{6X}{q\pi^2}\prod\limits_{p\mid \frac{q}{g}}\frac{1}{1-p^{-2}}+O_{\varepsilon}(X^{1/2}+q^{1/2+\varepsilon}),\label{squarefree error}
\end{align}
where to write the error term above from \cite[Equation (1.3)]{Nunes}, we are crudely using that $1\le g\le q$.

When $N$ is odd, using the Chinese remainder theorem, we get integers $a, \ b $ and $c$ such that the following equalities hold
\begin{align*}
   \sharp\, S_+(X,m,N)&= \sum\limits_{\substack{1\leq d\leq X\\ d\equiv m\pmod{N}\\ d\equiv 1\pmod{4}}}\mu^2(d)+ \sum\limits_{\substack{1\leq d\leq \frac{X}{4}\\ d\equiv m\pmod{N}\\ d\equiv 2\pmod{4}}}\mu^2(d)+ \sum\limits_{\substack{1\leq d\leq \frac{X}{4}\\ d\equiv m\pmod{N}\\ d\equiv 3\pmod{4}}}\mu^2(d)\nonumber\\
 &= \sum\limits_{\substack{1\leq d\leq X\\ d\equiv a\pmod{4N}}}\mu^2(d)+\sum\limits_{\substack{1\leq d\leq \frac{X}{4}\\ d\equiv b\pmod{4N}}}\mu^2(d)+\sum\limits_{\substack{1\leq d\leq \frac{X}{4}\\ n\equiv c\pmod{4N}}}\mu^2(d),
\end{align*}
where we are writing the first line above, simply by the definition of fundamental discriminants.

Note that $(a,4N)=(c,4N)=(m,N)$ and $(b,4N)=2(m,N)$. Moreover, the conditions on $m$ and $N$ in Lemma~\ref{Horie} ensures that every prime dividing $4N$ also divides $\frac{4N}{2(m,N)}$. Therefore, using $(\ref{squarefree error})$, the proof follows when $N$ is odd. Now let us consider the case when $N$ is even.

If $N\equiv 0\pmod{4}$ and $m\equiv1\pmod{4}$, then for any $d\equiv m\pmod{N}$ we have $d\equiv 1\pmod{4}$. Then, we can write
\begin{align*}
    \sharp\, S_+(X,m,N)&=\sum\limits_{\substack{1\leq d\leq X\\ d\equiv m\pmod{N}}}\mu^2(d).
\end{align*}

On the other hand, if $N\equiv 0\pmod{16}$ and $m\equiv 8\pmod{16}$ or $m\equiv12\pmod{16}$, then for any fundamental discriminant $d$ with $d\equiv m\pmod{N}$, we have $\frac{d}{4}\equiv2\pmod{4}$ or $\frac{d}{4}\equiv3\pmod{4}$ respectively. In this case, we can write
\begin{align*}
    \sharp\, S_+(X,m,N)&=\sum\limits_{\substack{1\leq d\leq \frac{X}{4}\\ d\equiv \frac{m}{4}\pmod{\frac{N}{4}}}}\mu^2(d),
\end{align*}
again, we are writing the equality above, by the definition of fundamental discriminants.

In both cases of even $N$, again, the proof follows from \eqref{squarefree error}, noting that the conditions on $m$ and $N$ in Lemma~\ref{Horie} imply that every odd prime dividing $N$ also divides $\frac{N}{(m,N)}$. 
\end{proof}

\subsubsection{Counting binary cubic forms with congruence conditions}
Let $M$ be the set of equivalence classes of binary integral cubic forms under the action of $\mathrm{GL}_2(\mathbb{Z})$. Note that the discriminant is invariant under the $\mathrm{GL}_2(\mathbb{Z})$ action. For any $F\in M$, we denote $\Delta(F)$ be the discriminant of $F$. Let $\Delta^+(X)$ (resp. $\Delta^-(X)$) denote the set of $F\in M$ with a positive (resp. negative) discriminant with $|\Delta(F)| \le X$.

Fix a set of exponents $\{\alpha_p:p~\mathrm{prime}\}$. For any prime $p$, we consider a subset $S_{p^{\alpha_p}}$ of $(\mathbb{Z}/p^{\alpha_p}\mathbb{Z})^4$. Let $E_{p^{\alpha_P}}$ be the set of $F\in M$, whose reduction modulo $p^{\alpha_p}$ is in $S_{p^{\alpha_p}}$. For any integer $n$, we then consider 
\begin{equation}\label{eqn:en}
E_n=\bigcap\limits_{\substack{p\mid n\\p~\mathrm{prime}}}E_{p^{\alpha_p}},\quad\mathrm{and}\quad E=\bigcap\limits_{p~\mathrm{prime}}E_{p^{\alpha_p}}.
\end{equation}

Au usual, we say that an integral cubic form $F\in M$ is in $E_n$, if $F \pmod n \in E_n$. Then, we want to count the number of $F\in M$ with the congruence conditions determined by $E_n$. To do that, certainly, we need to denote the \textit{local densities}, 
\begin{equation}\label{eq:sptp}
s(p)=\frac{\sharp\, S_{p^{\alpha_p}}}{p^{4\alpha_p}},\quad \mathrm{and} \quad t(p)=1-s(p).
\end{equation}

Let us now set $H^+=\frac{\pi^2}{72},\ H^-=\frac{\pi^2}{24}$, and for any real number $c>0$, denote the function 
$$E_c(X)=\exp\left(-c\sqrt{\log{X}\log{\log{X}}}\right).$$ 

To count the equivalence classes in $M$ with the congruence conditions determined by $E$, Belabas~\cite[Theorem 2.3]{Belabas} provided the following recipe.

\begin{lemma}\cite[Theorem 2.3]{Belabas}\label{belabas 2.3} 
Suppose that there exists some $\alpha>0$, so that $\alpha_p\leq 4\alpha$ for any prime $p$. Also, assume that there exists some $u>1$ such that the following holds;
\begin{itemize}
    \item[(i)] $t(p)= O(p^{-u}),$
    \item[(ii)] the forms in $E_{p^{\alpha_p}}$ are non zero modulo $p$,
    \item[(iii)] the number of classes in $\Delta^{\pm}(X)$ that are not in $E_p^{\alpha_p}$, is $O(Xp^{-u})$.
\end{itemize}
Then for any $0<c<c_0=\frac{u-1}{4\sqrt{\alpha^2+1}}$, the number of irreducible primitive classes of binary cubic forms in $\Delta^{\pm}(X)\cap E$ is given by
$$H^{\pm} X\prod_{p}s(p) +O_c(XE_c(X)).$$
\end{lemma}

Using the proof of Lemma~\ref{belabas 2.3} above, we now derive an explicit version of Lemma~\ref{Horie}, with particular attention to the dependence of the error term on the moduli $N$. Although our result is inspired by the argument in \cite[Theorem 2.3]{Belabas}, we cannot apply that theorem directly, as our moduli $N$ does not satisfy the necessary conditions stated there. Instead, we adapt the argument used in the proof of \cite[Theorem 2.3]{Belabas} to suit our setting.

\begin{proposition}\label{error}
 Let $m$ and $N$ be any integers as in Lemma~\ref{Horie}. Fix a real number $\frac{1}{4\sqrt{5}}-\frac{1}{16}<c_0<\frac{1}{4\sqrt{5}}$, and set 
$$\varepsilon=\lambda-c_0,\quad \mathrm{where}\quad \lambda=\frac{1}{4\sqrt{5}}.$$
Then we have the following estimate 
 \begin{equation}\label{eq:explicitNH}
   \sum\limits_{\substack{d\in S_+ (X,m,N)}}3^{r_3(d)}=\frac{4X}{\pi^2\varphi(N)}\prod\limits_{\substack{p\mid N\\ p \text{ prime}}}\frac{q}{p+1}+O\left(N^2XE_{c}(X)\right),
 \end{equation}
where $q$ is as in Lemma~\ref{S(X) error}, and $c=\frac{\lambda}{2}(1-16\varepsilon)>0$.
\end{proposition}

\subsection{Proof of Proposition~\ref{error}}
We denote $\alpha_p=\nu_p(N)$ for any prime $p\mid N$. Otherwise, we denote denote $\alpha_p=2$ if $p$ is odd, and $\alpha_p=4$ if $p=2$. Writing $N=\prod_{p\mid N} p^{\alpha_p}$ and using the Chinese remainder theorem, we find integers $m_p$, such that 
$$m\pmod{N}=m_p\pmod{p^{\alpha_p}},\quad \text{for all}\quad p\mid N.$$

Now, let $E_{p^{\alpha_p}}$ be the set of forms in $M$ whose discriminants are $\equiv m_p \pmod {p^{\alpha_p}}$, when $p\mid N$. Otherwise, set $E_{p^{\alpha_p}}$ be the forms in $M$ whose discriminants are $\not\equiv 0 \pmod {p^{\alpha_p}}$ when $p$ is odd, and are $\equiv 1,7,8,9,12$ or $13 \pmod {16}$, when $p=2$. Certainly, we then need to count the forms in the intersection of $E:= \bigcap\limits_{p\ \mathrm{prime}} E_{p^{\alpha_p}}$, and $\Delta^{\pm}(X)$. The same set-up is also considered in \cite[Page 22]{MR0958035}.

At this point, we argue exactly as done in the proof of Lemma~\ref{belabas 2.3}. To get the main term in \eqref{eq:explicitNH} (like the one appearing in \cite[Theorem 2.3]{Belabas}), we need to know all the local densities, as defined in \eqref{eq:sptp}. For our case, a description of them is readily available in \cite[Pages 22, 23]{MR0958035}. Alternatively, the main term can also be obtained directly from Lemma~\ref{Horie}.

The main task is to get the error term in \eqref{eq:explicitNH}. It is essentially sum of the error terms coming from Lemma~\ref{S(X) error} and the proof of \cite[Theorem 3.1]{Belabas}, with our choices of $u,\alpha$ and the constants appearing in $(\textit{i})$ and $(\textit{iii})$ of Lemma~\ref{belabas 2.3}. 

We go for the naive options, i.e, take $u=2$ along with the following constants in Lemma~\ref{belabas 2.3}:
\begin{equation}
C_1=2N^2 \quad \mathrm{in}\quad (\textit{i}),\quad \mathrm{and}\quad C_3=C N^2 \quad \mathrm{in}\quad (\textit{iii}),
\end{equation}
where $C>1$ is a sufficiently large constant, so that \cite[Proposition 1]{MR0491593} holds with constant $C$. 

Note that the choice of $C_1$ works, because we have $t(p)\le 2p^{-2}$ for any $p>N$ (following the proof of \cite[Theorem 2.3]{Belabas}). Also, the choice of $C_2$ works, because we can apply \cite[Proposition 1]{MR0491593} for primes $p>N$, and crudely use the estimate of $|\Delta^{\pm}(X)|$ (as mentioned on \cite[Page 4]{Belabas}) for primes $p\le N$. Moreover, the condition at $(ii)$ is satisfied, by the definitions of our sets $E_{p^{\alpha_p}}$.

Now, to apply Lemma~\ref{belabas 2.3} we first need $\alpha$. Note that $\alpha$ is crudely $\log N$ here, which would give us a poor saving in the error term. However, this is very much avoidable, as we shall see shortly below. We also need to figure out the dependence of $c$ (equivalently, $N$) in the constant appearing in the error term in Lemma~\ref{belabas 2.3}. To work these things out, let us revisit below the proof of Lemma~\ref{belabas 2.3}.
\subsubsection{On the sieving}
For a parameter $Y$ to be chosen suitably, consider the function $f=f^{+}$ as in \cite[Corollary 2.2]{Belabas}, and denote $P_Y$ to be the product of all primes $\le Y$. In our case, we have the following estimate for any $(r,P_Y)=1$, 
\begin{equation}\label{eq:fr}
f(r)= H^{+}\prod_{p\le Y} s(p)\prod_{p|r}t(p)X + O\left(X^{15/16+\varepsilon} N^{1/4}\prod_{p\mid rP_Y}p\right),
\end{equation}
and also, we can write the error term above as
\begin{equation}\label{eq:errorinfr}
O\left(X^{15/16+\varepsilon}\cdot N^{1/4}r \cdot \exp(2Y)\right).
\end{equation}

We now introduce another parameter \( Y < Z \), and invoke condition \((iii)\) to eliminate the contribution from the forms not lying in $E_{p^{\alpha_p}}$ for some primes $p > Z$, at the expense of only a negligible cost. To be more precise, we can write our desired quantity, i.e., \cite[Equation (3)]{Belabas} as \begin{equation}\label{eq:notinZ}
f(1) - \sum_{k > 1} \sum_{\substack{p_1 < \cdots < p_k \\ Y < p_i < Z}} (-1)^{k-1} f(p_1 \cdots p_k) + O(N^2 Z^{-1}).
\end{equation}

We now introduce another parameter
$K$, and write the main term in \eqref{eq:notinZ} as
\begin{equation}\label{eq:afterK}
f(1) - \sum_{1 \leq k < K} \sum_{\substack{p_1 < \cdots < p_k \\ Y < p_i < Z}} (-1)^{k-1} f(p_1 \cdots p_k)
+ O\left( \sum_{\substack{p_1 < \cdots < p_K \\ Y < p_i < Z}} f(p_1 \cdots p_K) \right).
\end{equation}

Now, we argue exactly as in \cite[Equation~(5)]{Belabas}, and use \eqref{eq:fr}, \eqref{eq:errorinfr} combined with $(i)$, to conclude that the term in \eqref{eq:afterK} is equals to
\begin{align*}
H^{+} \prod_{p \leq Y} s(p) X &\left( 1 + \sum_{k > 1}
\sum_{\substack{p_1 < \cdots < p_k \\ Y < p_i < Z}} (-1)^k t(p_1 \cdots p_k) \right)\\
&\quad+ O \left( X^{15/16 + \varepsilon}\cdot N^{1/4}K\cdot \exp(2Y)
\Big( \sum_{p < Z} p \Big)^K + X \Big( N^2 \sum_{p > Y} p^{-2}\Big)^K \right).
\end{align*}

Expanding and collecting the $t(p)$ terms in the main term above, we deduce that the term in \eqref{eq:notinZ} is
\[
H^\pm \prod_{p < Z} s(p) X + O \left(
X^{15/16 + \varepsilon}\cdot N^{1/4}K \cdot \exp(2Y) Z^{2K}
+ X \left( N^2 Y^{-1} \right)^K \right).
\]
Again, using $(i)$, we finally derive that \eqref{eq:notinZ} is
\[
H^\pm \prod s(p) X +O\left(
X^{15/16 + \varepsilon}\cdot N^{1/4} K \cdot \exp(2Y)  Z^{2K}
+ X \left( N^2 Y^{-1} \right)^K+N^2XZ^{-1} \right).
\]

Dividing both sides by $X$ to make things simpler, the task is now to optimize the following term
\begin{equation}\label{eq:finalerror}
N^2 \left(
X^{-1/16 + \varepsilon}\cdot K \cdot \exp(2Y) Z^{2K}
+ Z^{-1} \right)+(N^2 Y^{-1})^K,
\end{equation}
for a suitable choice of parameters $Y, Z $ and $K$.

\subsubsection{Concluding the proof}
We go with the same choices of $Y,Z$ and $K$ as in \cite[Pages 5, 6]{Belabas}, i.e., we choose
$$Y = \frac{\log X}{\log \log \log X}, \quad K=\lambda(1-16\varepsilon) \left( \frac{\log X}{\log \log X}\right)^{1/2},\quad \mathrm{and}$$
$$Z=\exp\left(\lambda (\log X \log \log X)^{1/2}\right).$$
Then, the first term in \eqref{eq:finalerror} is immediately dominated by $O\left(N^2XE_{c_0}(X)\right)$; because it is essentially at most $N^2$ times the error term appearing in \cite[Equation (4)]{Belabas}.

On the other hand, note that $c<c_0$, hence the estimate in \eqref{eq:explicitNH} is trivial if 
$$N^2\ge \exp\left(c_0 \sqrt{\log X \log \log X}\right).$$

If $N$ does not satisfy the above, then the second term $(N^2 Y^{-1})^K$ in \eqref{eq:finalerror} is clearly dominated by $\exp\left(-c\sqrt{\log X \log \log X}\right)$, concluding the proof.
\qed

\section{Indivisibility over real quadratic fields}

Let us recall that we denote $S$ to be the set of all square-free numbers. As a consequence of the results established in the section earlier, we have the following two estimates, which will play a key role in the proof of both Theorem~\ref{indizuka} and Theorem~\ref{squrefree}.

\begin{lemma}\label{prop}
We have,
   $$ \sum\limits_{\substack{1\leq d\leq X\\d\in S}}3^{r_3(d)}=\frac{8}{\pi^2}X+O\left(XE_c(X)\right),\quad \mathrm{for~some}\quad c>0.$$
\end{lemma}

\begin{proof}
We can write
    \begin{align*}
\sum\limits_{\substack{1\leq d\leq X\\d\in S}}3^{r_3(d)}&=\sum\limits_{\substack{1\leq d\leq X\\ d\in S \\ d\equiv 1\pmod{4}}}3^{r_3(d)}+\sum\limits_{\substack{1\leq d\leq X\\ d\in S \\ d\not\equiv 1\pmod{4}}}3^{r_3(d)}\\
&=\sum\limits_{\substack{1\leq d\leq X\\ d\in S_+ \\ d\equiv 1\pmod{4}}}3^{r_3(d)}+\sum\limits_{\substack{1\leq d\leq 4X\\ d\in S_+ \\ d\not\equiv 1\pmod{4}}}3^{r_3(d)}\\
&=\sum\limits_{d\in S_+(X,1,4)}3^{r_3(d)}+\sum\limits_{d\in S_+(4X,1,1)}3^{r_3(d)}-\sum\limits_{d\in S_+(4X,1,4)}3^{r_3(d)},
\end{align*}
where we deduce the second line from the first, applying \eqref{eq:simplefact} on the right-most term. The proof now follows from Proposition~\ref{error}.
\end{proof}

Consequently, we deduce the following estimate.
\begin{lemma}\label{avergesum in N}
    We have,
$$\sum\limits_{\substack{1\leq n\leq X}}3^{r_3(n)}=\frac{4}{3}X+o(X)$$
\end{lemma}
\begin{proof}
We fix a parameter $0<\varepsilon<1/2$ to be chosen suitably. Then, writing any natural number $n$ as $dm^2$ with $d\in S$ and $m\in \mathbb{N}$, we can split the following sum as 
\begin{equation}\label{eq:decomp}
\sum_{1 \leq n \leq X} 3^{r_3(n)} = \sum_{\substack{dm^2 \leq X \\ d \in S}} 3^{r_3(dm^2)} 
= \sum_{\substack{d \leq X \\ d \in S}} \left\lfloor \sqrt{\frac{X}{d}} \right\rfloor 3^{r_3(d)}=S_1+S_2,
\end{equation}
where the sums $S_1$ and $S_2$ are defined by
$$S_1= \sum_{1 \leq k \leq X^{\frac{1}{2} - \varepsilon}} 
k \sum_{\substack{d \in S \\ \frac{X}{(k+1)^2} < d \leq \frac{X}{k^2}}} 3^{r_3(d)},\quad \mathrm{and}$$
$$S_2=\sum_{X^{\frac{1}{2} - \varepsilon} \leq k \leq X^{\frac{1}{2}}} 
k \sum_{\substack{d \in S \\ \frac{X}{(k+1)^2} < d \leq \frac{X}{k^2}}} 3^{r_3(d)}.$$

By Lemma~\ref{prop}, we have
\begin{align*}
S_1&= \sum_{1 \leq k \leq X^{\frac{1}{2} - \varepsilon}} 
    k~\Big( 
    \sum_{\substack{d \in S \\ 1 < d \leq \frac{X}{k^2}}} 3^{r_3(d)} 
    - \sum_{\substack{d \in S \\ 1 < d \leq \frac{X}{(k+1)^2}}} 3^{r_3(d)} 
    \Big)\\
    &=\frac{8}{\pi^2} X \sum_{1 \leq k \leq X^{\frac{1}{2} - \varepsilon}} 
     k \left( 
    \frac{1}{k^2} - \frac{1}{(k+1)^2} 
    \right) + O\left(X \log X E_c(X^{2\varepsilon})\right),\quad \mathrm{for~some}\quad c>0,
\end{align*}
where we are writing the error term above, noting that we have $X^{2\varepsilon}\le X/k^2$, in the range of $k$ above.

On the other hand, we estimate $S_2$ rather crudely as 
\begin{align*}
    S_2 
    &\leq X^{\frac{1}{2}} 
    \sum_{X^{\frac{1}{2} - \varepsilon} \leq k \leq X^{\frac{1}{2}}}~\sum_{\substack{d \in S \\ \frac{X}{(k+1)^2} < d \leq \frac{X}{k^2}}} 
    3^{r_3(d)} \leq X^{\frac{1}{2}} 
    \sum_{\substack{d \in S \\ d \leq X^{2\varepsilon}}} 
    3^{r_3(d)} =O\left( X^{\frac{1}{2} + 2\varepsilon}\right),
\end{align*}
where again, we use Lemma~\ref{prop} to get the last estimate above. 

The proof concludes, noting that
\begin{equation}\label{eq:kindentity}
    X\sum_{1 \leq k \leq X^{\frac{1}{2} - \varepsilon}} 
     k \left( 
    \frac{1}{k^2} - \frac{1}{(k+1)^2} 
    \right)=\frac{\pi^2}{6}X+o(X),
\end{equation}
and by plugging the estimates of $S_1$ and $S_2$ in \eqref{eq:decomp}, taking any $0<\varepsilon<1/4$.
\end{proof}
\subsection{Proof of Theorem \ref{indizuka}}
Let $k$ be any natural number. Note that
$$\sum\limits_{\substack{1\leq n\leq X \\ r_3(n)<k}}1+3^k\left(X-\sum\limits_{\substack{1\leq n\leq X \\  r_3(n)<k}}1\right)\leq\sum\limits_{1\leq n\leq X}3^{r_3(n)}.$$

Therefore, Lemma~\ref{avergesum in N} immediately shows that the set 
$\left\{n\in\mathbb{N}: 3^k\nmid h_{\Q(\sqrt{n})}\right\}$ has a density at least 
\begin{equation}\label{eqn:density}
\frac{1}{3^{k}-1}\left(3^k-\frac{4}{3}\right)=1-\frac{1}{3(3^k-1)},
\end{equation}
in the set of natural numbers. Let us now consider the set 
$$\mathbb{A}=\left\{n\in\mathbb{N}: 3^k\nmid h_{\Q(\sqrt{n})}, \  3^k\nmid h_{\Q(\sqrt{n+1})},\ldots, \ 3^k\nmid h_{\Q(\sqrt{n+m})}\right\},$$
and denote $c=1-\frac{m+1}{3(3^k-1)}$.

If $\mathbb{A}$ has density strictly less than \( c \) in $\mathbb{N}$, then \( \mathbb{N} \setminus \mathbb{A} \) has a density strictly more than \( 1-c \) in $\mathbb{N}$. Note that, for each \( n \in \mathbb{N} \setminus \mathbb{A} \), there exists some $0\leq i\leq m$ such that the class number of \( \mathbb{Q}(\sqrt{n+i}) \) is divisible by \( 3^k \). 

Hence, the set $\left\{n \in \mathbb{N} : 3^k \mid h_{\mathbb{Q}(\sqrt{n})}\right\}$ clearly has a density strictly more than \( \frac{1-c}{m+1}:= \frac{1}{3(3^k-1)}\) in $\mathbb{N}$, which is a contradiction to \eqref{eqn:density}, and the proof concludes.
\qed

\subsection{Proof of Theorem~\ref{squrefree}}
The proof follows by arguing similarly to the proof of Theorem~\ref{indizuka}. The only difference is that we use Lemma~\ref{prop}, instead of Lemma~\ref{avergesum in N}.
\qed

\section{Over imaginary quadratic fields}
In this section, we prove analogous results for imaginary quadratic fields. Let us first note the following result of Scholz \cite{MR1581309}, which allows us the transition from real to imaginary quadratic fields.

\begin{theorem}{(Scholz \cite{MR1581309})}\label{scholz}
    Let $n$ be any natural number. If we denote
$$r=r_3(n),\quad \mathrm{and}\quad s=r_3(-3n).$$
Then, we have $r\leq s\leq r+1$.
\end{theorem}
\begin{proposition}\label{prop:alldivby3}
           For any natural number $k$, the set 
$$\left\{n\in\mathbb{N}: n\equiv 0\pmod{3}, \  3^k\nmid h_{\Q(\sqrt{n})}\right\},$$ 
has a density at least $1-\frac{1}{3(3^k-1)}$ in $3\mathbb{N}$.
    \end{proposition}
    \begin{proof}
We can split the following sum as
\begin{equation}\label{eqn:decomp}
 \sum_{\substack{1 \leq n \leq X \\ n \equiv 0 \pmod{3}}} 3^{r_3(n)}= \sum_{\substack{1 \leq dt^2 \leq X \\ dt^2 \equiv 0 \pmod{3} \\ d \in S}} 3^{r_3(dt^2)}=S_1+S_2-S_3,
 \end{equation}
where $S_1$ denotes the contributions from $3\mid d,~S_2$ be the contribution from $3\mid t$, and $S_3$ be the contribution from $3\mid (d,t)$. Now, the task is to estimate $S_1,S_2$ and $S_3$.

\subsubsection*{Estimation of $S_1$} We proceed as in the proof of Lemma~\ref{avergesum in N}, and write 
\begin{align}
S_1&=\sum\limits_{\substack{1\leq d\leq X\\ 3\mid d, \ d\in S}}\left\lfloor\sqrt{\frac{X}{d}}\right\rfloor 3^{r_3(d)}=\sum\limits_{1\leq k\leq  X^{\frac{1}{2}-\varepsilon}}k\sum\limits_{\substack{d\in S, \ 3\mid d\\\frac{X}{(k+1)^2}<d\leq \frac{X}{k^2}}}3^{r_3(d)}\label{3divd'}\\
&\qquad\qquad\qquad\qquad\qquad\qquad\qquad+\sum\limits_{ X^{\frac{1}{2}-\varepsilon}\leq k\leq X^{\frac{1}{2}}}k\sum\limits_{\substack{d\in S, \ 3\mid d\\\frac{X}{(k+1)^2}<d\leq \frac{X}{k^2}}}3^{r_3(d)},\label{3divd}
    \end{align}
for a parameter $0<\varepsilon<1/2$, which we shall choose suitably.
    
To estimate the last sum in \eqref{3divd'}, let us first write
\begin{equation}\label{eqn:3div}
\sum\limits_{\substack{1\leq d\leq X\\ d\in S, \ 3\mid d}}3^{r_3(d)}=\sum\limits_{\substack{1\leq d\leq X\\ d\in S, \ 3\mid d\\ d\equiv 1\pmod{4}}}3^{r_3(d)}+\sum\limits_{\substack{1\leq d\leq X\\ d\in S, \ 3\mid d\\ d\not\equiv1\pmod{4}}}3^{r_3(d)}.
\end{equation}

Let us note that
\begin{align*}
    \sum\limits_{\substack{1\leq d\leq X\\ d\in S, \ 3\mid d\\ d\equiv 1\pmod{4}}}3^{r_3(d)}&=\sum\limits_{\substack{1\leq d\leq X\\ d\in S_+\\ d\equiv 1\pmod{4}}}3^{r_3(d)}-\sum\limits_{\substack{1\leq d\leq X\\ d\in S_+, \ d\equiv1\pmod{3}\\ d\equiv 1\pmod{4}}}3^{r_3(d)}\\
    &\hspace{4cm}-\sum\limits_{\substack{1\leq d\leq X\\ d\in S_+, \ d\equiv2\pmod{3}\\ d\equiv 1\pmod{4}}}3^{r_3(d)}.
\end{align*}

In particular, applying Proposition~\ref{error}, the first sum in the right-hand side of (\ref{eqn:3div}) is
\begin{equation}\label{eq:7firsteq}
   \sum\limits_{\substack{1\leq d\leq X\\ d\in S, \ 3\mid d\\ d\equiv 1\pmod{4}}}3^{r_3(d)}=\frac{2}{3\pi^2}X+O(XE_c(X)),\quad \mathrm{for~some}\quad c>0.
\end{equation}

On the other hand, one can write the last sum on the right-hand side of (\ref{eqn:3div}) as
\begin{align}
    &\sum_{\substack{1 \leq d \leq X \\ d \in S, \ 3 \mid d \\ d \not\equiv 1 \pmod{4}}} 3^{r_3(d)}= \sum_{\substack{1 \leq d \leq 4X \\ d \in S_+, \ 3 \mid d \\ d \not\equiv 1 \pmod{4}}} 3^{r_3(d)} = \sum_{\substack{1 \leq d \leq 4X \\ d \in S_+, \ 3 \mid d}} 3^{r_3(d)} 
       - \sum_{\substack{1 \leq d \leq 4X \\ d \in S_+, \ 3 \mid d \\ d \equiv 1 \pmod{4}}} 3^{r_3(d)}\label{eq:identityagain} \\
    &\qquad\qquad= \sum_{\substack{1 \leq d \leq 4X \\ d \in S_+}} 3^{r_3(d)} 
       - \sum_{\substack{1 \leq d \leq 4X \\ d \in S_+, \ d \equiv 1 \pmod{3}}} 3^{r_3(d)} -\sum_{\substack{1 \leq d \leq 4X \\ d \in S_+, \ d \equiv 2 \pmod{3}}} 3^{r_3(d)}\nonumber \\
    &\qquad \qquad\qquad\qquad\qquad\qquad\qquad \qquad\qquad- \sum_{\substack{1 \leq d \leq 4X \\ d \in S_+, \ 3 \mid d \\ d \equiv 1 \pmod{4}}} 3^{r_3(d)}\label{eq:lastsum},
\end{align}
where we are again using the simple fact from \eqref{eq:simplefact} to get the first equality in \eqref{eq:identityagain}. 

Again applying Proposition~\ref{error}, and estimating the sum at \eqref{eq:lastsum} by \eqref{eq:7firsteq}, we derive
\begin{align*}
\sum\limits_{\substack{1 \leq d \leq X \\ d \in S, \ 3 \mid d \\ d \not\equiv 1 \pmod{4}}} 3^{r_3(d)}=\frac{4}{3\pi^2}X+O(XE_c(X)),\quad \mathrm{for~some}\quad c>0.
\end{align*}

In particular, the sum at (\ref{eqn:3div}) is now given by
$$\sum\limits_{\substack{1\leq d\leq X\\  d\in S, \ 3\mid d }}3^{r_3(d)}=\frac{2}{\pi^2}X+O(XE_c(X)).$$ 

Now we estimate the sums on the right-hand side of (\ref{3divd}) as in the proof of Lemma~\ref{avergesum in N}, and again plugging \eqref{eq:kindentity}, we derive
\begin{align}
   S_1=\frac{2}{\pi^2}\cdot \frac{\pi^2}{6}X+o(X).
    \label{3divd first term}
\end{align}

\subsubsection*{Estimation of $S_2$}
Again, as in the estimation of $S_1$, or as in the proof of Lemma~\ref{avergesum in N}, we can write
\begin{align*}
    S_2&=\sum\limits_{\substack{1\leq 9dk^2\leq X\\ d\in S}}3^{r_3(d)}
    =\sum\limits_{\substack{1\leq d\leq X\\ d\in S}}\left\lfloor\sqrt{\frac{X}{9d}}\right\rfloor 3^{r_3(d)}\\
    &=\sum\limits_{1\leq k\leq  X^{\frac{1}{2}-\varepsilon}}k\sum\limits_{\substack{d\in S \\\frac{X}{9(k+1)^2}<d\leq \frac{X}{9k^2}}}3^{r_3(d)}+\sum\limits_{ X^{\frac{1}{2}-\varepsilon}\leq k\leq X^{\frac{1}{2}}}k\sum\limits_{\substack{d\in S\\\frac{X}{9(k+1)^2}<d\leq \frac{X}{9k^2}}}3^{r_3(d)}.
\end{align*}

At this point, we again argue as in the proof of Lemma~\ref{avergesum in N}, and derive
\begin{align}
     S_2&=\frac{8}{9\pi^2}\frac{\pi^2}{6}X+o(X).\label{3divd second term}
\end{align}

\subsubsection*{Estimation of $S_3$}
We can write
 \begin{align*}
   S_3&=\sum\limits_{\substack{1\leq 9dk^2\leq X\\ k\in\mathbb{N},\ 3\mid d, \ d\in S}} 3^{r_3(d)}=\sum\limits_{\substack{1\leq d\leq X\\3\mid d, \ d\in S}}\left\lfloor\sqrt{\frac{X}{9d}}\right\rfloor3^{r_3(d)}\\
        &=\sum\limits_{1\leq k\leq  X^{\frac{1}{2}-\varepsilon}}k\sum\limits_{\substack{d\in S, \ 3\mid d\\\frac{X}{9(k+1)^2}<d\leq \frac{X}{9k^2}}}3^{r_3(d)}+\sum\limits_{ X^{\frac{1}{2}-\varepsilon}\leq k\leq X^{\frac{1}{2}}}k\sum\limits_{\substack{d\in S, \ 3\mid d\\\frac{X}{9(k+1)^2}<d\leq \frac{X}{9k^2}}}3^{r_3(d)}.
\end{align*}

Once again, we argue as in the proof of Lemma~\ref{avergesum in N}, while now using the estimate from \eqref{3divd first term}. We derive
\begin{align}
  S_3=\frac{2}{9\pi^2}\frac{\pi^2}{6}X+o(X).\label{3divd third term}  
\end{align}

Finally, substituting \eqref{3divd first term}, \eqref{3divd second term} and \eqref{3divd third term} in \eqref{eqn:decomp}, we deduce
\begin{align*}
    \lim\limits_{X\to\infty} \ \frac{\sum\limits_{\substack{1\leq n\leq X\\ 3\mid n}}3^{r_3(n)}}{\sharp\, \left\{1\leq n\leq X : 3\mid n\right\}} = \frac{4}{3}.
\end{align*}

From this point, the proof follows by the same argument as we used Lemma~\ref{avergesum in N} to derive \eqref{eqn:density}.
\end{proof}
\subsection{Proof of Theorem~\ref{indizukaimaginary}}
Let us consider the set 
$$\mathbb{A}=\left\{n\in 3\mathbb{N}: 3^k\nmid h_{\Q(\sqrt{n})}, \  3^k\nmid h_{\Q(\sqrt{n+3})},\ldots, \ 3^k\nmid h_{\Q(\sqrt{n+3m})}\right\}.$$

By the same argument as in the proof of Theorem~\ref{indizuka}, note that Proposition~\ref{prop:alldivby3} implies that, $\mathbb{A}$ has density at least $1-\frac{m+1}{3(3^k-1)}$ in the set $3\mathbb{N}$. The proof finally concludes, applying Theorem~\ref{scholz}.
\qed

\section{Over bi-quadratic fields}\label{sec:biquad} 
In this section, we prove Theorem~\ref{biquad}. Let us first recall the necessary results.
\subsection{Preliminaries}
Let $x$ be a positive real number. Denote $D(x)$ be the density of the set of natural numbers $n$ such that $\frac{n}{\varphi(n)}\leq x$. 

\begin{lemma}[Erd\"os~\cite{MR16078}]\label{erdos}
    Let $0<\varepsilon<1$ be a fixed real number. Then, for any sufficiently large $x$, we have
    \begin{align*}
        \exp(-\exp((1+\varepsilon)ax))<1-D(x)< \exp(-\exp((1-\varepsilon)ax)),\quad a=e^{-\gamma},  
      \end{align*}
      where $\gamma$ is Euler's constant.
\end{lemma}

Cornell~\cite{MR1123163} proved the following relation between the class number of an imaginary bi-quadratic field and its quadratic subfields. 

\begin{lemma}{\cite[Lemma 2]{MR1123163}}\label{cornell}
    The class number of an imaginary bi-quadratic field is
either the product of the class numbers of the three quadratic subfields or half of that number. 
\end{lemma}

Byeon~\cite[Proposition~3.1]{MR2073286} proved the existence of certain pairs of quadratic fields whose class numbers are not divisible by $3$. This serves as a key tool for the proof of Theorem~\ref{biquad}. It is formulated for a fixed \(t\). 
To establish Theorem~\ref{biquad}, it will be necessary to account for the dependence on \(t\). To this end, we present the following refinement of \cite[Proposition~3.1]{MR2073286}.




\begin{lemma}\label{modified byeon}
   Let $t$ be any odd negative square-free integer, and $m, N$ be two natural numbers as in Lemma~\ref{Horie}. Assume that $N$ is odd and $(m,t)=1$. Then, we have 
\begin{align*}
      \sharp\, &\left\{d\in S_+(Y,m,tN):3\nmid h_{\Q(\sqrt{d})}, \ 3\nmid h_{\Q(\sqrt{td})}\right\}\\
&\qquad\qquad\qquad\quad \geq\frac{1}{3} \sharp\, S_+(Y,m,tN)-C_0 N^2t^2YE_{c}(Y),
\end{align*}
for some constant $C_0>0$, independent to $t,m$ and $N$, where $c$ is as in Proposition~\ref{error}.
\end{lemma}
\begin{proof}

Note that our imposed conditions ensure that whenever \(m\) and \(N\) satisfy the congruence conditions in Lemma~\ref{Horie}, the numbers \(m\) and $-tN$, as well as $-mt$ and \(t^2N\), also satisfy them. 
We then proceed exactly as in the proof of \cite[Proposition~3.1]{MR2073286}, and instead of applying \cite[Lemma~2.2]{MR2073286}, we apply its explicit version due to Proposition~\ref{error}. 
This completes the proof.
    \end{proof}

    

\subsection{Proof of Theorem~\ref{biquad}} It follows from Lemma \ref{densS+} and Lemma~\ref{S(X) error} (taking $m=1$ and $N=4$) that, the density of negative odd fundamental discriminants whose class numbers are not divisible by $3$, is at least $\frac{1}{2}\left(\frac{2}{\pi^2}\right)=\frac{1}{\pi^2}$ in the negative integers. Let us denote $G$ be the set of such negative odd fundamental discriminants. 

For any $m$ and odd negative square-free $t\in G$ with $(m,3t)=1$, using Lemma~\ref{S(X) error} (with $\varepsilon=1/2$) and Lemma~\ref{modified byeon} (for $N=3$), we derive
\begin{align*}
  &\sharp\, \left\{d\in S_+(Y,m,3t):3\nmid h_{\Q(\sqrt{d})}, \ 3\nmid h_{\Q(\sqrt{td})}\right\}\\
  &\hspace{2.5cm}\geq\frac{Y}{2\pi^2\varphi(-t)}\prod\limits_{p\mid 3t}\frac{p}{p+1}-C t^2YE_{c}(Y),   
\end{align*}
for some constant $C>0$, independent to $m$ and $t$.

In particular,
\begin{align*}
    &\sum_{\substack{1 \leq m < -3t \\ (m, 3t) = 1}} 
    \sharp\, \left\{ d \in S_+(Y, m, 3t) : 3 \nmid h_{\mathbb{Q}(\sqrt{d})}, \ 3 \nmid h_{\mathbb{Q}(\sqrt{td})} \right\} \\
    &\hspace{2cm} \geq \frac{Y}{\pi^2} \prod_{p \mid 3t} \frac{p}{p+1} -C\varphi(-3t)t^2YE_{c}(Y) \\
    &\hspace{2cm} \geq \frac{3Y}{4\pi^2} \frac{\varphi(-t)}{-t}-C\varphi(-3t)t^2YE_{c}(Y).
\end{align*}

Let $\varepsilon=1-\frac{\ln{\ln{\frac{\pi^2}{\pi^2-9}}}}{a}, \ a=e^{-\gamma}$, and $x_0$ be some real number satisfying the following
  \begin{align*}
  1-D(x_0)<e^{-e^{(1-\varepsilon)ax_{0}}}<e^{-e^{(1-\varepsilon)a}}=\frac{\pi^2-9}{\pi^2}.    
  \end{align*}
  
  Lemma~\ref{erdos} shows that $G\cap T(x_0)$ has a positive density in $-\mathbb{N}$, where $T(x_0)=\{n\in-\mathbb{N}:\frac{\varphi(-n)}{-n}\geq \frac{1}{x_0}\}$. In particular, we can find some $c_{x_0}>0$ such that the following holds for any sufficiently large $X$; 
  $$\sharp\, \left\{t\in G\cap T(x_0): 1\leq -t\leq X\right\}\geq c_{x_0} X.$$ 
  
Recall the function $H(X,Y)$ from \eqref{eq:H}. By Lemma~\ref{cornell}, for any sufficiently large $X$ and $Y$, we now have
\begin{align*}
    H (X,Y)&\geq \sum\limits_{\substack{t\in G\cap T(x_0)\\ 1\leq -t\leq X}} \left(\frac{3Y}{4\pi^2}\frac{\varphi(-t)}{(-t)}+3Ct^3YE_{c}(Y)\right)\\
    &\geq  \frac{3Y}{4\pi^2}\frac{c_{x_0}}{x_0}X-C'X^4YE_{c}(Y),
\end{align*}
for some constant $C'>0$. The proof concludes, as the imposed condition on $Y$ ensures that the second term above is $o(XY)$.
\qed
\section{Further questions and remarks}\label{sec:remarks}
For any interval $[x,y]$, let us denote 
$$D[x,y] = \left\{ n \in [x,y] \cap \mathbb{N} \right\} \setminus \left\{ n^2 : n \in \mathbb{N} \right\},$$ 
and
$$S[x,y] = \left\{ n \in D[x,y] : 3 \nmid h_{\Q(\sqrt{n+i})}, \ \forall i \in \{0,1,2,3,4\} \right\}.
$$ 

Using SageMath, we generated the following table, presenting a statistic for the sizes of the sets $D[x,y]$ and $S[x,y]$ over various ranges. This is to illustrate Theorem~\ref{indizuka}, at least for $k=1$ and $m=4$.
\begin{table}[h!]      
\centering
\begin{minipage}[b]{0.65\textwidth}
\centering

\vspace{0.2cm}
\begin{tabular}{|c|c|c|}
            \hline
            \( x \text{ to } y  \) & \( \sharp\, D[x,y] \)  & $ \sharp\, S[x,y]$  \\
            \hline
            2 \text{ to } 2000 & 1851 & 1392 \\
            2001 \text{ to } 4000 & 1839 & 1358 \\
            4001 \text{ to} 6000 & 1831 & 1312\\
             6001 \text{ to } 8000 & 1819 & 1274\\
            8001 \text{ to } 10000 & 1819 & 1256 \\
          10001 \text{ to } 12000 & 1803 & 1213\\
            12001 \text{ to } 14000 & 1814 & 1279\\
              14001 \text{ to } 16000 & 1801 & 1190  \\
               16001 \text{ to } 18000 & 1810 & 1231 \\
               18001 \text{ to } 20000 & 1809 & 1266\\
                 20001 \text{ to } 22000 & 1796 & 1196\\
                     22001 \text{ to } 24000 & 1807 & 1230 \\
                      24001 \text{ to } 26000 & 1799 & 1203 \\
                        26001 \text{ to } 28000 & 1790 & 1186 \\
                        28001 \text{ to } 30000 & 1801 & 1218 \\
            \hline
    \end{tabular}\label{tab:corollary_p_k}
    \end{minipage}
    \end{table}

\subsection{Some further natural questions}
Let \( m \) be a odd natural number. Consider the sum  
\begin{equation}\label{eqn:comp}
\sum_{0 < D \leq X} \prod_{p \mid m} p^{r_p(D)},
\end{equation}
where the product is over all prime divisors \( p \) of \( m \). A study of this sum helps to understand the simultaneous indivisibility by the prime divisors of \( m \). To get an upper bound for (\ref{eqn:comp}), we perhaps require a generalization of \cite[Conjecture 4]{BM19}, for the density of quadratic fields with simultaneously prescribed $p$ and $q$-ranks, for any two distinct primes $p$ and $q$.

Theorem~\ref{indizukaimaginary} is certainly valid for an indivisibility by $3^{k}$ with $k\geq 2$. A natural question would be the case $k=1$?

In Lemma~\ref{Horie}, we have the average $3$-torsion over $S_+(X,m,N)$. When $3$ is replaced by any other prime $p$ and Conjecture~\ref{cohen} holds, then a result analogous to Theorem~\ref{indizuka} holds with $p^k$. It would then give a positive answer to \cite[Question~1]{CS22}, for primes $p>3$ and many suitable $t$.

\section*{Acknowledgements}
We would like to thank IISER-TVM for providing excellent working conditions. We would like to thank Jayanta Manoharmayum, Pasupulati Sunil Kumar and Jishu Das for their helpful comments. Muneeswaran wishes to thank CSIR for the funding. SK's research was supported by SERB grant CRG/2023/009035. SB is supported by an institute fellowship from IISER-TVM and ARC grant DP230100530, with appreciation to UNSW Sydney for providing an excellent work environment.

\end{document}